\documentclass[a4paper,12pt]{article}

\usepackage[english]{babel}
\usepackage[utf8]{inputenc}
\usepackage[T1]{fontenc}

\usepackage{amsmath,amssymb,color}
\usepackage{ntheorem}

\usepackage{hyperref}

\usepackage[alphabetic]{amsrefs}

\usepackage{graphicx}

\usepackage{verbatim}
\usepackage{multirow}
\usepackage{float}

\usepackage{geometry}
\theoremstyle{plain}
\theoremheaderfont{\normalfont\bfseries}
\theorembodyfont{\itshape}
\theoremseparator{}
\theoremnumbering{arabic}
\theoremsymbol{}
\newtheorem{defn}{Definition}[section]

\newcommand{\p}{\partial}

\title{{\it Quasi}-Random Discrete Ordinates Method for Transport Problems}
\author{$^a$Pedro H. A. Konzen, $^b$Leonardo F. Guidi, $^c$Thomas Richter\\
  \small{$^{a,b}$ IME, UFRGS, Porto Alegre, Brazil; $^c$ IAN, OVGU, Magdeburg, Germany}\\
  \small{$^a$\texttt{pedro.konzen@ufrgs.br}, $^b$\texttt{guidi@mat.ufrgs.br}, $^c$\texttt{thomas.richter@ovgu.de}}x}
\date{}
\begin{document}

\maketitle

\section*{Abstract}

The {\it quasi}-random discrete ordinates method (QRDOM) is here proposed for the approximation of transport problems. Its central idea is to explore a {\it quasi} Monte Carlo integration within the classical source iteration technique. It preserves the main characteristics of the discrete ordinates method, but it has the advantage of providing mitigated ray effect solutions. The QRDOM is discussed in details for applications to one-group transport problems with isotropic scattering in rectangular domains. The method is tested against benchmark problems for which DOM solutions are known to suffer from the ray effects. The numerical experiments indicate that the QRDOM provides accurate results and it demands less discrete ordinates per source iteration when compared against the classical DOM.

\tableofcontents

\section{Introduction}

In this paper the following one-group transport problem in an
isotropic medium with reflective boundary conditions is considered~\cite{Lewis1984a}:
\begin{subequations}\label{eq:tp}
  \begin{align}
    \forall \Omega\in S^2:~&&\Omega\cdot\nabla \psi(\pmb{x},\Omega) + \sigma_t\psi(\pmb{x},\Omega) &= \frac{\sigma_s}{4\pi}\int_{S^2} \psi(\pmb{x},\Omega'')\,d\Omega'' + Q(\pmb{x}), &\forall \pmb{x}&\in\mathcal{D},\label{eq:te}\\
    \forall \Omega\in S^2, \pmb{n}\cdot\Omega < 0:~&&\psi(\pmb{x},\Omega) &= \rho \psi(\pmb{x},\Omega') + Q_b(\pmb{x}),& \pmb{x}&\in\Gamma,\label{eq:tbc}    
  \end{align}
\end{subequations}
where $S^2 := \{(\mu, \eta, \xi):~\mu^2 + \eta^2 + \xi^2 = 1\}$ is the sphere in $\mathbb{R}^3$, $\nabla$ is the gradient operator in $\mathbb{R}^3$, $\psi(\pmb{x},\Omega)$ is the intensity at point $\pmb{x}$ in the domain $\mathcal{D}\subset\mathbb{R}^3$ in the direction $\Omega\in S^2$, $\sigma_t$ and $\sigma_s$ are, respectively, the total and scattering macroscopic cross sections, $Q(\pmb{x})$ and $Q_b(\pmb{x})$ are, respectively, the sources in $\mathcal{D}$ and on its boundary $\Gamma$, $\pmb{n}$ is the unit outer normal on $\Gamma$, $\Omega' := \Omega - 2(\pmb{n}\cdot\Omega)\pmb{n}$ is the reflected direction of $\Omega$ on $\Gamma$, and $\rho$ is the reflective coefficient.

One of the most widely used techniques to solve \eqref{eq:tp} is the
classical discrete ordinates method (DOM) (see, for instance,
\cite[Ch. 3]{Lewis1984a}), which consists in approximating the
integral term in the left-hand side of the equation \eqref{eq:te} by
using an appropriate quadrature set $\{\Omega_i, w_i\}_{i=1}^M$. This
leads to the approximation of the integro-differential equation by a
system of partial differential equations on the discrete ordinates
$\{\psi(\pmb{x},\Omega_i)\}_{i=1}^M$, which can be solved by a variety
of classical discretization methods and easily be integrated in CFD codes. It is well known that the quality of the DOM solution depends on the choice of the quadrature set \cite{Koch2004a,Hunter2013a}. In particular, for transport problems with discontinuities in the source, with discontinuities in the cross sections, or on non-convex geometries, the DOM approximation may produce unrealistic oscillatory solutions known as the ray effects \cite{Chai1993a,Morel2003a}.

Ray effects can be mitigated by increasing the number of discrete
ordinates \cite{Li2003a}, at the expense of additional computational
costs. In fact, in order to compute approximations with a large number
of discrete ordinates, one needs a robust computational
implementation, otherwise it may not be feasible due to a large memory demand. Many remedies for ray effects have been proposed in the last decades. Integral methods to the transport problem \cite{Loyalka1975a,Altac2004a,Azevedo2018a} and the Modified DOM Method \cite{Ramankutty1997a} are known to produce accurate results, but they are not as straightforward to integrate to Computational Fluid Dynamic codes as the DOM scheme. It is also known that the choice of the DOM quadrature set plays an important role in the accuracy of the solution \cite{AbuShumays2001a,Barichello2016a}. Alternatively, Adaptive Discrete Ordinates schemes have been proposed \cite{Stone2007a, Jarrel2010a}. Recently, the Frame Rotation Method (FRM) \cite{Tencer2016a} has been proposed. Given a quadrature set, the FRM computes the transport problem solution as the simple mean of DOM solutions obtained from random rotations of the quadrature set on the sphere.

In this paper the {\it quasi}-random discrete ordinates method (QRDOM) is proposed for the approximation of the transport problem solution with mitigated ray effects. Its central idea is to explore a {\it quasi} Monte Carlo integration \cite{Leobacher2014} within the classical source iteration technique. It admits a parallelizable computational implementation enhanced by a convergence acceleration. Although it is here introduced as an alternative technique to solve \eqref{eq:tp}, it may be adapted to more general cases of multigroup transport problems in anisotropic mediums. The major advantage of the proposed QRDOM is the mitigation of the ray effects without the loss of the good characteristics of the DOM.

In Section \ref{sec:QRDOM} the fundamentals of the QRDOM are presented. In Section \ref{sec:QRDOM_rect_domains} the QRDOM application to transport problems in rectangular domains is discussed in details. Then in Section \ref{sec:results} selected numerical experiments with the QRDOM applied to known benchmark problems are presented. Finally, in Section \ref{sec:final} final considerations are given.

\section{{\it Quasi}-Random Discrete Ordinates Method}\label{sec:QRDOM}

The proposed {\it Quasi}-Random Discrete Ordinates Method (QRDOM) is based on the idea of the well-known {\it Quasi}-Monte Carlo Method for integration. More explicitly, we assume that problem \eqref{eq:tp} can be approximated by
\begin{subequations}\label{eq:mctp}
  \begin{align}
    1\leq i \leq M:~&&\Omega_i\cdot\nabla \psi_i + \sigma_t\psi_i &= \frac{\sigma_s}{M}\sum_{k=1}^M \psi_k + Q(\pmb{x}),&\forall \pmb{x}&\in\mathcal{D},\label{eq:mcte}\\
    \forall \pmb{n}\cdot\Omega_i < 0:~&&\psi(\pmb{x},\Omega_i) &= \rho \psi(\pmb{x},\Omega_i') + Q_b(\pmb{x}), ~&\pmb{x}&\in\Gamma,\label{eq:mctbc}
  \end{align}
\end{subequations}
where $\{\Omega_i\}_{i=1}^{M}\subset S^2$ is a given {\it quasi}-random finite sequence of discrete directions, and $\psi_i = \psi(\pmb{x})$ is the approximation of $\psi(\pmb{x},\Omega_i)$.

This approximation of the transport problem~\eqref{eq:tp} poses the issue to select {\it a priori} the total number of discrete ordinates $M$. It is expected that $M$ must be of order of thousands for the most problems. Then, the issue of solving~\eqref{eq:mctp} poses a computational challenge.

With this in mind, we propose the following numerical iterative strategy to approximate the solution of the transport problem. Firstly, let us denote the scalar flux by
\begin{equation}
  \Psi(\pmb{x}) := \frac{1}{4\pi}\int_{S^2} \psi(\pmb{x},\Omega)\,d\Omega,
\end{equation}
and assume $\Psi^{(1)}$ is a given initial approximation of it at the so-called \emph{epoch one}. Let us also denote by $qr:\mathbb{N}\to S^2$ a {\it quasi}-random (low discrepancy) generator of discrete directions $qr(i) = \Omega_i\in S^2$. Then, we compute the $\Psi^{(l)}$ approximation of the scalar flux at the $l$-th epoch by solving
\begin{subequations}\label{eq:qrtp}
  \begin{align}
    && \Omega_i\cdot\nabla \psi_i^{(l)} + \sigma_t\psi_i^{(l)} &=
    \sigma_s\Psi^{(l-1)} + Q, &&\text{in }\mathcal{D},\\
    \text{if}~\pmb{n}\cdot\Omega_i < 0:~&&\psi^{(l)}(\Omega_i) &= \rho
    \psi^{(l)}(\Omega_i') + Q_b, &&\text{on }\Gamma,
  \end{align}
\end{subequations}  
for each $i=m^{(l-1)},m^{(l-1)}+1,m^{(l-1)}+2,\dotsc,m^{(l)}-1$, and by accumulatively computing 
\begin{equation}
  \Psi^{(l)} = \frac{1}{M^{(l)}}\sum_{k=m^{(l-1)}}^{m^{(l)}-1}\psi_k^{(l)},
\end{equation}
where $m^{(1)}=1$, $m^{(l)}-1$ is the lower index of the {\it quasi}-random sequence greater than $m^{(l-1)}$ for which $\Psi^{(l)}$ converges to a given precision at this epoch, and $M^{(l)} = m^{(l)}-m^{(l-1)}$. This kind of source iteration procedure continues until a desired convergence is achieved.

We should note that this proposed iterative procedure has the advantage that problem~\eqref{eq:qrtp} involves just $qr(i) = \Omega_i$ and its reflected $\Omega_i'$ directions on the boundary. Moreover, by this strategy the total number of discrete ordinates $M^{(l)}$ at epoch $l$ is determined {\it a posteriori} by a chosen tolerance.

In the next section, we restrict ourselves to the transport problem in rectangular domains and, in this context, we address the issues of building the {\it quasi}-random sequence of discrete direction and approximating the solution of~\eqref{eq:qrtp}.

\section{Applications in rectangular domains}\label{sec:QRDOM_rect_domains}

\begin{figure}[h]
  \centering
  \includegraphics[width=0.7\textwidth]{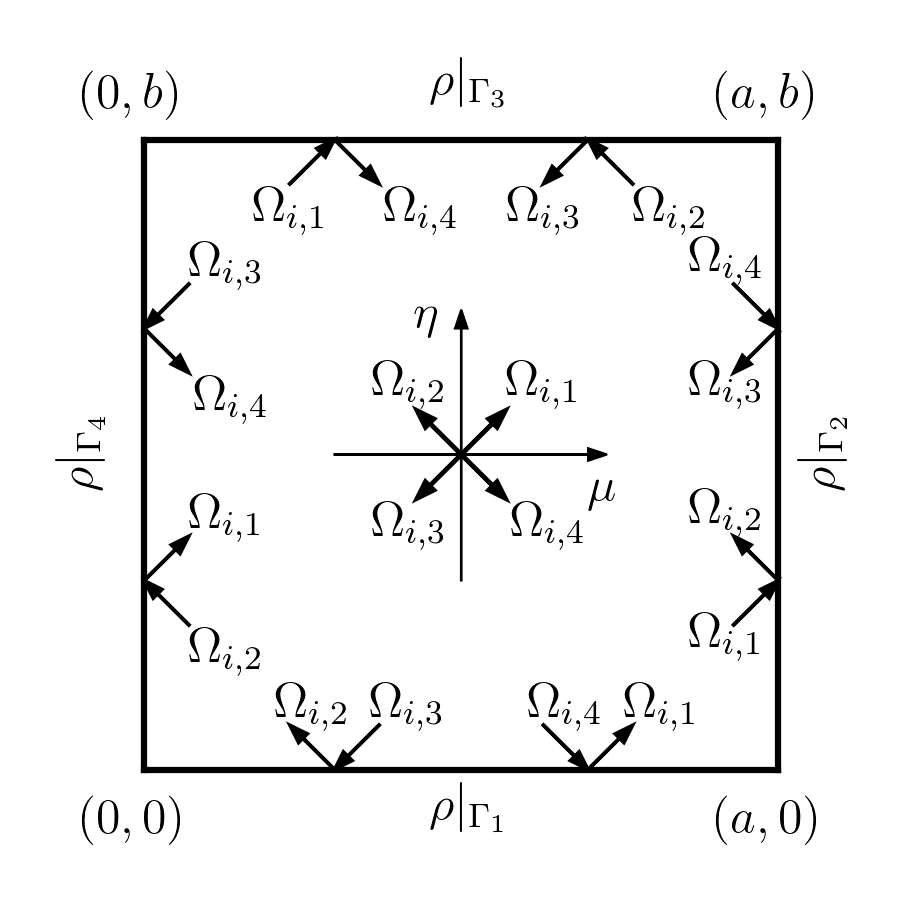}
  \caption{Illustration of the rectangular with incident and reflected directions on the boundary.}
  \label{fig:domain}
\end{figure}

Here, we explore the proposed QRDOM method for applications in
rectangular domains, more explicitly, we assume $\mathcal{D}=(0,
a)\times (0, b)$ and also assume symmetry in the $z$-coordinate. In
this case, the complexity of problem~\eqref{eq:qrtp} can be further
reduced as follows. Let us denote by $S^2_1$ the first octant of the sphere, i.e. $S^2_1 := \{(\mu,\eta,\xi):~\mu^2+\eta^2+\xi^2=1, \mu>0, \eta>0, \xi>0\}$. Then, let us also consider $qr_1:\mathbb{N}\to S^2_1$ a {\it quasi}-random generator of discrete directions on $S^2_1$, i.e. $qr_1(i) = \Omega_{i,1} =  (\mu_{i,1}, \eta_{i,1}, \xi_{i,1})\in S^2_1$. Associated with $qr_1(i)$ we also consider its boundary reflected directions: $\Omega_{i,2} = (-\mu_{i,2}, \eta_{i,2}, \xi_{i,2})$, $\Omega_{i,3} = (-\mu_{i,3}, -\eta_{i,3}, \xi_{i,3})$, and $\Omega_{i,4} = (\mu_{i,4}, -\eta_{i,4}, \xi_{i,4})$ (see Figure \ref{fig:domain}).

Therefore, the QRDOM procedure presented in the previous section can be slightly adapted as follows: the $\Psi^{(l)}$ approximation at epoch $l$ is built from the solution to
\begin{subequations}\label{eq:qrtp_rect}
  \begin{align}
    j=1,2,3,4:~&\mu_{i,j}\frac{\p \psi_{i,j}^{(l)}}{\p x} + \eta_{i,j}\frac{\p \psi_{i,j}^{(l)}}{\p y} + \sigma_t\psi_{i,j}^{(l)} = \sigma_s\Psi^{(l-1)} + Q(\pmb{x}), \forall \pmb{x}\in\mathcal{D}, \label{eq:qrtp_rect_eq}\\
    &\psi_{i,4}^{(l)} = \rho \psi_{i,1}^{(l)} + Q_b,\qquad \psi_{i,3}^{(l)} = \rho \psi_{i,2}^{(l)} + Q_b, \text{on}~\Gamma_1:=[0,a]\times \{0\}, \label{eq:qrtp_rect_bc1}\\
    &\psi_{i,1}^{(l)} = \rho \psi_{i,2}^{(l)} + Q_b, \qquad \psi_{i,4}^{(l)} = \rho \psi_{i,3}^{(l)} + Q_b, \text{on}~\Gamma_2:=\{a\}\times [0, b],  \label{eq:qrtp_rect_bc2}\\
    &\psi_{i,1}^{(l)} = \rho \psi_{i,4}^{(l)} + Q_b, \qquad \psi_{i,2}^{(l)} = \rho \psi_{i,3}^{(l)} + Q_b, \text{on}~\Gamma_3:=[0,a]\times \{b\},  \label{eq:qrtp_rect_bc3}\\
    &\psi_{i,3}^{(l)} = \rho \psi_{i,4}^{(l)} + Q_b, \qquad \psi_{i,2}^{(l)} = \rho \psi_{i,1}^{(l)} + Q_b, \text{on}~\Gamma_4:=\{0\}\times [0,b],  \label{eq:qrtp_rect_bc4}
  \end{align}
\end{subequations}  
for each $i=m^{(l-1)},m^{(l-1)}+1,m^{(l-1)}+2,\dotsc,m^{(l)}-1$, $qr_1(i) = \Omega_{i,1}$, and by accumulatively computing
\begin{equation}
  \Psi^{(l)} = \frac{1}{M^{(l)}}\sum_{k=m^{(l-1)}}^{m^{(l)}-1}\Psi^{(l)}_k,
\end{equation}
where $\Psi^{(l)}_k := (\psi_{k,1}^{(l)}+\psi_{k,2}^{(l)}+\psi_{k,3}^{(l)}+\psi_{k,4}^{(l)})/4$ is the scalar flux of the $k$-th sample, and again $m^{(l)}-1$ is the lower index of the {\it quasi}-random sequence greater than $m^{(l-1)}$ for which $\Psi^{(l)}$ converges to a given precision at this epoch. 

In the following Subsection \ref{subsec:qr}  $qr_1$ is presented and,
based hereon, a convergence acceleration is proposed in Subsection \ref{subsec:conv_accel} together with the definition of a convergence criterion.

\subsection{{\it Quasi}-random generator}\label{subsec:qr}

The choice of the {\it quasi}-random generator is a key point for the QRDOM, since it will directly impact the convergence of the iterative procedure. Here, $qr_1$ is built from the bi-dimensional $\mathcal{S}_{2,3}$ {\it quasi}-random reverse Halton sequence~\cite{Vandewoestyne2006a}, which is denoted as $rh:\mathbb{N}\to (0,1)\times (0,1)$, $rh(i)=(rh_1(i), rh_2(i))$. More explicitly, $qr_1:\mathbb{N}\to S_1^2$, $qr(i) = \Omega_i = (\mu_i, \eta_i, \xi_i)$ is taken as
\begin{equation}\label{eq:param-S2}
  \begin{split}
    \mu_i &= \sin(\arccos(1-rh_1(i)))\cos(rh_2(i)\pi/2),\\
    \eta_i &= \sin(\arccos(1-rh_1(i)))\sin(rh_2(i)\pi/2),\\
    \xi_i &= \cos(\arccos(1-rh_1(i))).
  \end{split}
\end{equation}

The discrepancy is a measure of how far, in a certain sense, a finite
sequence of elements in $[0,1)^d$ is from a uniformly distributed modulo one sequence in the same region. A mathematically precise description for the intuitive notion of uniformly distributed points is provided for infinite sequences of them contained in finite intervals. According to this intuitive notion, given an infinite sequence of elements in $[0,1)^d$, the ratio of those contained in some region is proportional to the region's size.

The definition of uniform distribution is constructed  by means of semi-open  intervals like $[\boldsymbol{a},\boldsymbol{b})=[a_1,b_1)\times\ldots\times[a_d,b_d)$, where indices refer to components of $\boldsymbol{a}$ and $ \boldsymbol{b}$, both elements of $[0,1)^d$.  Let $\mathcal{S}=\left(\boldsymbol{x}^{(1)},\boldsymbol{x}^{(2)},\ldots\right)$ be an infinite sequence of elements in $[0,1)^d$, $\left.\mathcal{S}\right|_N$ the finite sequence given by the first $N$ elements of $\mathcal{S}$ and $\lambda_d\left([\boldsymbol{a},\boldsymbol{b})\right)$ the $d$-dimensional Lebesgue measure of some interval $[\boldsymbol{a},\boldsymbol{b})\subseteq[0,1)^d$. In order to define the concept of uniform distribution for sequences, it is necessary to know the cardinality of the set  of the indices  of the elements of $\mathcal{S}$ that  belongs to $[\boldsymbol{a},\boldsymbol{b})$. It will be symbolized by  $A\left([\boldsymbol{a},\boldsymbol{b}),\mathcal{S}\right)$,
\begin{equation}
A\left([\boldsymbol{a},\boldsymbol{b}),\mathcal{S}\right):=
\left|\left\{n\in\mathbb{N} : \boldsymbol{x_n}\in[\boldsymbol{a},\boldsymbol{b})\right\}\right|.
\end{equation} 

\begin{defn} An infinite sequence $\mathcal{S}$ in $[0,1)^d$ is \emph{uniformly distributed modulo one}, if 
	\begin{equation}
	\lim_{N\rightarrow\infty}\frac{A\left([\boldsymbol{a},\boldsymbol{b}),\left.\mathcal{S}\right|_N\right)}{N}=\lambda_{d}\left([\boldsymbol{a},\boldsymbol{b})\right)
	\end{equation}
	for every interval of the form $[\boldsymbol{a},\boldsymbol{b})\subseteq[0,1)^d$.
\end{defn}

It is a well know result that given an infinite sequence $\mathcal{S}=\left(\boldsymbol{x}^{(1)},\boldsymbol{x}^{(2)},\ldots\right)$ of elements in  $[0,1)^d$, and any Riemann integrable function $f:[0,1]^d\rightarrow\mathbb{R}$, the equality
\begin{equation}
\lim_{N\rightarrow\infty}\frac{1}{N}\sum_{n=1}^{N}f\left(\boldsymbol{x}^{(n)}\right)=\int_{[0,1]^d}f(\boldsymbol{\nu})d\boldsymbol{\nu}\label{sumint}
\end{equation} holds if and only if $\mathcal{S}$ is an infinite uniformly distributed modulo one sequence. 

The notion of discrepancy gives rise to some definitions of discrepancy (see \cite{Leobacher2014} and references therein), hereforth  we will consider the \emph{star discrepancy}.

\begin{defn} Let $\mathcal{S}$ be a sequence in $[0,1)^d$. The \emph{star discrepancy} of this set, $\mathcal{D}^{*}_{N}\left(\mathcal{S}\right)$  is defined as the number 
	\begin{equation}
	\mathcal{D}^{*}_{N}\left(\mathcal{S}\right):=\sup_{\boldsymbol{x}_\in[0,1)^d}\left|\frac{A\left([\boldsymbol{0},\boldsymbol{x}),\left.\mathcal{S}\right|_N\right)}{N}-\lambda_{d}\left([\boldsymbol{0},\boldsymbol{x})\right)\right|.
	\end{equation}
\end{defn}
The star discrepancy allows a form of error estimate in the
approximation of integrals by partial sums in the limit (\ref{sumint})
where the contribution of the size of the finite sequence is
independent of the integrand. Now, let us restrict ourselves to
bi-dimensional sequences. Let $f:[0,1)^{2}\rightarrow\mathbb{C}$ be a
  function whose mixed partial derivatives of order one are all
  continuous and let 
$\left\Vert f\right\Vert_{2,1}$ be a norm given by
\begin{equation}
\left\Vert f\right\Vert_{2,1}:=\sum_{s\subseteq\{1,2\}}\int_{[0,1]^{|s|}}\left|\frac{\partial^{|s|}f}{\partial \boldsymbol{\nu}_s}\left(\boldsymbol{\nu}_s,\boldsymbol{1}\right)\right|d\boldsymbol{\nu}^s
\end{equation} where  $\left(\boldsymbol{\nu}_s,\boldsymbol{1}\right)=\left(\iota_1,\iota_2\right)$  is the \emph{anchored} argument in subset $s$:
\begin{equation}
\iota_i:=\left\{\begin{array}{ll} \nu_i & \text{if }i\in s\\ 1 & \text{if }i\notin s \end{array}\right..
\end{equation} This choice of norm allows the following version of Hlawka-Koksma inequality \cite{Leobacher2014} for the error when a finite sequence $\left.\mathcal{S}\right|_N$ of knots in $[0,1)^2$ is used to approximate the integral of $f$
\begin{equation}
\left|\frac{1}{N}\sum_{n=1}^{N}f\left(\boldsymbol{x}^{(n)}\right)-\int_{[0,1]^2}f(\boldsymbol{\nu})d\boldsymbol{\nu}\right|\leq\left\Vert f\right\Vert_{2,1}\mathcal{D}^{*}_{N}\left(\mathcal{S}\right).\label{HK-ineq}
\end{equation} Hence, the speed of convergence strongly depends on the discrepancy's rate of decrease for sequence $\left.\mathcal{S}\right|_N$. The trick is to chose a sequence whose discrepancy decays rapidly. A result from Schmidt \cite{Schmidt1972} states that any finite sequence $\left.\mathcal{S}\right|_N\in[0,1)^2$  has a lower bound  for its star discrepancy given by
\begin{equation}
\mathcal{D}^{*}_{N}\left(\mathcal{S}\right)\geq C\,\frac{\log N}{N},\label{Dstarlowerb}
\end{equation} where $C>0$ is a constant.
It is a known result that an infinite sequence  $\mathcal{S}\in [0,1)^{d}$ is uniformly distributed modulo one if and only if  $\lim_{N\rightarrow\infty}\mathcal{D}^{*}_{N}\left(\mathcal{S}\right)=0$ (see \cite{Leobacher2014}).
The lower bound (\ref{Dstarlowerb}) limits the rate of convergence of partial sum approximations (for $d=2$). Ideally, one should choose a sequence whose upper bound has a behavior as near as possible of that exhibited by the right-hand side of (\ref{Dstarlowerb}). The generalized Halton sequences are a common choice. If $\mathcal{S}$ is a generalized Halton sequence, its star discrepancy is bounded as
\begin{equation}
N\mathcal{D}^{*}_{N}(\mathcal{S})\leq c_{d}\log^{d}N + O\left(\log^{d-1}N\right),\label{Dstarupperb}
\end{equation}where $c_d$ is a constant which depends on the choice of a base of $d$ pairwise coprime numbers  and the particular permutation used to construct the sequence (see \cite{Vandewoestyne2006a}, \cite{Foure2009} and \cite{Leobacher2014}).

\begin{figure}[h]
  \centering
\includegraphics[width=0.7\textwidth]{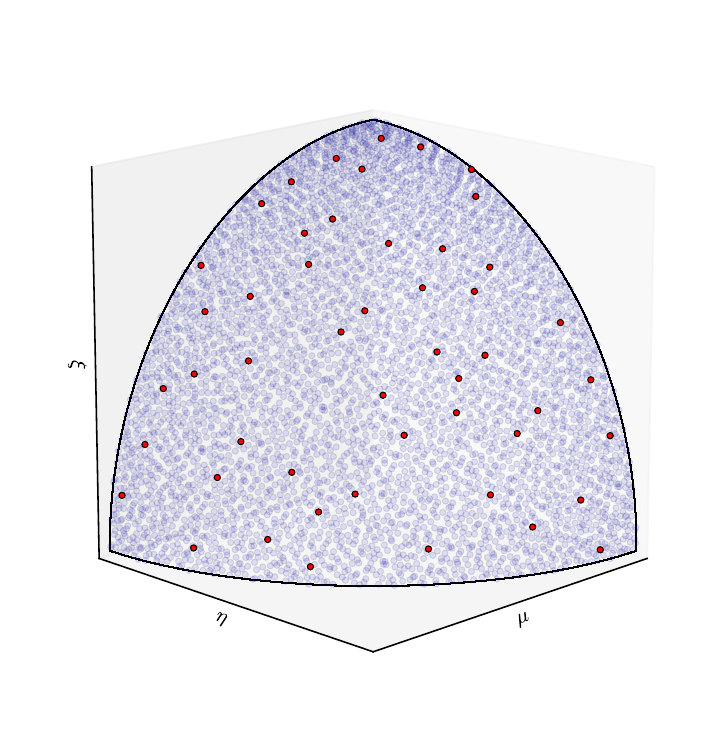}
  \caption{The first $5000$ {\it quasi}-random ordinate directions (samples) given by the $qr$ generator. The first $50$ samples are plotted as red circles.}
  \label{fig:samples}
\end{figure}

The parametrization of the first spherical octant given by (\ref{eq:param-S2}) in terms of $rh$ preserves uniform distribution. Therefore, as $rh$ produces a low discrepancy sequence in $(0,1)\times(0,1)$, so does  $qr$ in $S_1^2$. Figure \ref{fig:samples} shows the first $5000$ {\it quasi}-random ordinate directions given by the $qr$ generator.

\subsection{Convergence acceleration}\label{subsec:conv_accel}

The behavior of a goal functional of the scalar flux can be used to track the convergence of the QRDOM. More specifically, let $F^{(l)}_k = F(\Psi^{(l)}_k)$ denote the value of a given functional $F$ applied to the flux $\Psi^{(l)}_k$ of the $k$-th sample in the $l$-th epoch. Assuming further that $F$ is linear, it follows
\begin{equation}
  F(\Psi^{(l)}) = \frac{1}{M^{(l)}}\sum_{k=m^{(l-1)}}^{m^{(l)}-1}F_k^{(l)}.\label{linear_funct_F}
\end{equation}
Here one needs to determine the $m^{(l)}-1$ index. This index is defined from the convergence behavior of the preliminary $F^{(p)}\left(\Psi^{(l)}\right)$, \begin{equation}
F^{(p)}\left(\Psi^{(l)}\right)=\frac{1}{p}\sum_{k=m^{(l-1)}}^{m^{(l-1)}-1+p}F_k^{(l)},~p\geq 1.
\end{equation} 
Formally, the limit $p\rightarrow\infty$ gives the value of the
functional $F^{\infty}(\Psi)$ which is the value of the functional $F$
applied to the scalar flux of the exact solution of (\ref{eq:mcte}),
(\ref{eq:mctbc}). In light of this, the following strategy to define
$m^{(l)}$ is proposed. 

From Hlawka-Koksma inequality (\ref{HK-ineq}) and the bounds given by
(\ref{Dstarlowerb}) and (\ref{Dstarupperb}), the linear model function $\Phi$ can be considered for the behavior of $F^{(n)}\left(\Psi^{(l)}\right)$  with the number of samples $n$ taken from the $l$-th epoch,
\begin{equation}\label{eq:fun_to_fit}
\Phi(n;\left\{\gamma_{0},\gamma_{1},\gamma_{2}\right\})=\gamma_{0}+\gamma_{1}\dfrac{\log n}{n}+\gamma_{2}\dfrac{\log^{2}n}{n}.
\end{equation}
The estimate  $\tilde{F}^{\infty}\left(\Psi^{(l)}\right)$  for $F^{\infty}(\Psi)$ is given by  the weighted least square fit of $\Phi$ to the data given by $ \left\{\left(n,F^{(n)}\left(\Psi^{(l)}\right)\right)\right\}_{n=1}^{p^{*}}$ with weights\footnote{The choice of values for weights of the samples with $n<8$ is due to the non monotone behavior of the weight function in that region. } $w_{n}$,
\begin{equation}
w_{n}=\left\{\begin{array}{ll} 2 & \text{if }n<8,\\ \dfrac{n^{2}}{\log^{4}n} & \text{if }n\geq8. \end{array}\right.
\end{equation} 
More specifically, $\tilde{F}^{\infty}\left(\Psi^{(l)}\right)=\gamma_{0}$ and $p^{*}$ is the smallest number of samples that satisfies the convergence criterion 
\begin{equation}
\frac{\left|\tilde{F}^{\infty}\left(\Psi^{(l)}\right) - \tilde{F}^{\infty}\left(\Psi^{(l-1)}\right)\right|}{\left|\tilde{F}^{\infty}\left(\Psi^{(l)}\right)\right|} < \text{TOL}.\label{eq:tol_criterion}
\end{equation}
Then, the $l$-th epoch last index $m^{(l)}-1$ is given by 
\begin{equation}
  m^{(l)}=m^{(l-1)}+p^{*}.
\end{equation}

By imposing a linear goal functional, the QRDOM allows an efficient and trivial parallel implementation. This is the only reason of the linearity assumption. For any trial functional $G=G(\Psi^{(l)})$, the QRDOM solution is taken as its extrapolated estimate $\tilde{G}^\infty(\Psi^{(l)})$, which can be computed in a postprocessing step. Therefore and for the sake of simplicity, the tilde and the $\infty$ notation will be dropped out from now on.

\subsection{Implementation details}

The implementation of the proposed QRDOM method requires the application of a numerical method to solve problem \eqref{eq:qrtp_rect}. Here, the standard finite element method (FEM) with the streamline diffusion stabilization (SUPG) is applied. Briefly, the standard weak formulation of problem \eqref{eq:qrtp_rect} can be written as: for each $i$ find $u_i\in V := (H^1(\mathcal{D}))^4$ such that
\begin{equation}
  a(u_i,\varphi) + b_{\Gamma_-}(u_i,\varphi) = l(\varphi),\quad \forall\varphi\in V,
\end{equation}
where $u_i = (\psi_{i,1}, \psi_{i,2}, \psi_{i,3}, \psi_{i,4})$, $a(u_i,\varphi)$ and $l(\varphi)$ are the bilinear and the linear forms of the weak formulation equation \eqref{eq:qrtp_rect_eq}, respectively, and $b_{\Gamma_-}(u,\varphi)$ accounts for the weak form of the inflow boundary conditions \eqref{eq:qrtp_rect_bc1}-\eqref{eq:qrtp_rect_bc4}. Then, the discrete finite element problem reads: for each $i$ find $u_{i,h}\in V_h$ such that:
\begin{equation}\label{eq:qrdom_fem}
  a(u_{i,h},\varphi_h+\delta\mathcal{T}\varphi_h) + b_{\Gamma_-}(u_{i,h},\varphi_h+\delta\mathcal{T}\varphi_h) = l(\varphi_h+\delta\mathcal{T}\varphi_h),\quad \forall\varphi_h\in V_h,
\end{equation}
where $V_h\subset V$ is the finite element space of quadratic elements $Q_1$ build of a regular mesh, $\mathcal{T}\varphi_h := \Omega_i\cdot\nabla \phi_h$, and $\delta \approx 0.5h$ is the stabilization parameter, with $h$ denoting the element dimension.

One should expect that the QRDOM will demand the solution of thousands
of equations of type \eqref{eq:qrdom_fem} in each epoch until convergence is achieved with a common accepted tolerance. Therefore, a standard sequential algorithm may demand an excessive computational time. Alternatively, an implementation in a parallel MPI (master-slaves) paradigma has been developed. In the implemented code, the master processor instance controls the {\it quasi}-random generator, distributes the discrete directions to the slaves, and monitors the convergence. The slave instances receive the discrete directions from the master instance, solve the associated finite element problem \eqref{eq:qrdom_fem} and send the solution back to the master. As soon as the master receives a solution from a given slave, it sends to it a new discrete direction to be computed. Synchronization between all processor instances are required at each convergence checks. However it is not adequate to check the convergence at each new sample given the randomness of the procedure. In our numerical experiments (see Section \ref{sec:results}) checking the convergence at each thousand samples has been sufficient.

The computer code has been written in C++ with the help of the finite element toolkit Gascoigne 3D \cite{Gascoigne}. Moreover, the implemented code uses the GNU Science Library (GSL) \cite{GSL} for the reversed Halton {\it quasi}-random generator and the weighed least-square fitting.

\section{Numerical experiments}\label{sec:results}

In this section the performance of the QRDOM is discussed based on its application to benchmark problems. The selected problems share the rectangular computational domain $\mathcal{D} = (0,a)\times (0,b)$ as illustrated in Figure~\ref{fig:domain}.

\subsection{Benchmark problem 1: black walls}

The benchmark problem 1 has the parameters $a=b=2.5$, $\sigma_t\equiv 1$, $\sigma_s = 0.7-0.3(x^2/a^2+y^2/b^2)$, domain source $Q\equiv 0$, black walls $\rho_1=\rho_2=\rho_3=\rho_4=0$, and boundary sources $Q_{b1}\equiv 1$, $Q_{b2}=Q_{b3}=Q_{b4}\equiv 0$. It is well known that DOM solutions of this problem suffer of ray effects (see, for instance, \cite{Fiveland1984a,Li2003a,Truelove1987a, Ramankutty1997a,Tencer2016a} for studies on similar benchmark problems).

\begin{figure}[h]
  \centering
  \includegraphics[width=0.7\textwidth]{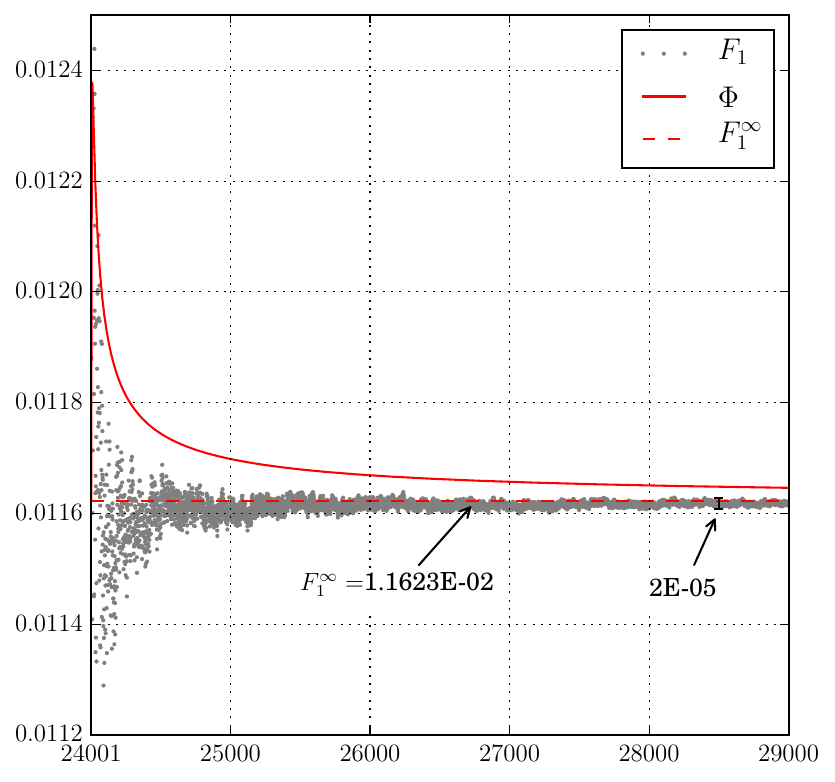}
  \caption{The goal functional values $F_1$ (black points) computed in the last epoch of the QRDOM applied to the benchmark problem 1. The fitted function $\Phi$ is plotted in solid red line, and the estimated value $F_1^{\infty}$ is plotted in dashed red line.}
  \label{fig:altac-conv}
\end{figure}

The application of the QRDOM to this problem was performed by assuming the total line scalar flux at $y=2.5$
\begin{equation}
  F_1(\Psi) := \int_{\Gamma_3} \Psi\,ds
\end{equation}
as the goal functional and with a tolerance of $TOL = 10^{-5}$ in the
convergence criterion
\eqref{eq:tol_criterion}. Table~\ref{tab:qrdom-altac-L3F} present the
computed $F_1$ and the scalar flux at the points $x=0$, $x=1.25$ and
$x=2.5$ in the top line $y=2.5$. Taking as a reference the values
reported in \cite[Table 4]{Altac2004a} one can confirm the QRDOM
convergence precision of at least $10^{-5}$. Figure
\ref{fig:altac-conv} shows the scatter plot of the $F_1(\Psi_k^{(l)})$
values in the last epoch of the QRDOM with the mesh of $256\times 256$
cells. The red solid line is the fitted function $\Phi(n)$ given in
\eqref{eq:fun_to_fit} and the red dashed line indicates its
extrapolation $F_1^\infty=1.1623\cdot 10^{-2}$. The error bar shows the minimum and the maximum values of $F_1(\Psi_k^{(k)})$ of the last 1000 samples.

\begin{table}[h]
  \centering
  \caption{Scalar flux values on the top wall computed with $QRDOM$ for the benchmark problem 1: $a=b=2.5$, $\sigma_s=0.7-0.3(x^2/a^2+y^2/b^2)$. Reference values~\cite[Table 4]{Altac2004a}: $\Psi(0,2.5) = 0.00891$, $\Psi(1.25,2.5)=0.01342$, $\Psi(2.5,2.5)=0.00758$. Total line scalar flux at $y=2.5$ as goal functional ($F_1$), $TOL = 10^{-4}$.}
  \begin{tabular}[H]{l|cccc}\hline
    \#cells & $F_1$ & $\Psi(0,2.5)$ & $\Psi(1.25,2.5)$ & $\Psi(2.5,2.5)$ \\\hline
    \multirow{1}{*}{$16\times 16$}   & $1.1629\cdot 10^{-2}$ & $8.9198\cdot 10^{-3}$ & $1.3447\cdot 10^{-2}$ & $7.5946\cdot 10^{-3}$ \\
    \multirow{1}{*}{$32\times 32$}   & $1.1626\cdot 10^{-2}$ & $8.9133\cdot 10^{-3}$ & $1.3433\cdot 10^{-2}$ & $7.5869\cdot 10^{-3}$ \\
    \multirow{1}{*}{$64\times 64$}   & $1.1624\cdot 10^{-2}$ & $8.9097\cdot 10^{-3}$ & $1.3430\cdot 10^{-2}$ & $7.5843\cdot 10^{-3}$ \\
    \multirow{1}{*}{$128\times 128$} & $1.1623\cdot 10^{-2}$ & $8.9069\cdot 10^{-3}$ & $1.3428\cdot 10^{-2}$ & $7.5813\cdot 10^{-3}$ \\
    \multirow{1}{*}{$256\times 256$} & $1.1623\cdot 10^{-2}$ & $8.9066\cdot 10^{-3}$ & $1.3427\cdot 10^{-2}$ & $7.5812\cdot 10^{-3}$ \\\hline
  \end{tabular}
  \label{tab:qrdom-altac-L3F}
\end{table}

\begin{figure}[h]
  \centering
  \includegraphics[width=0.7\textwidth]{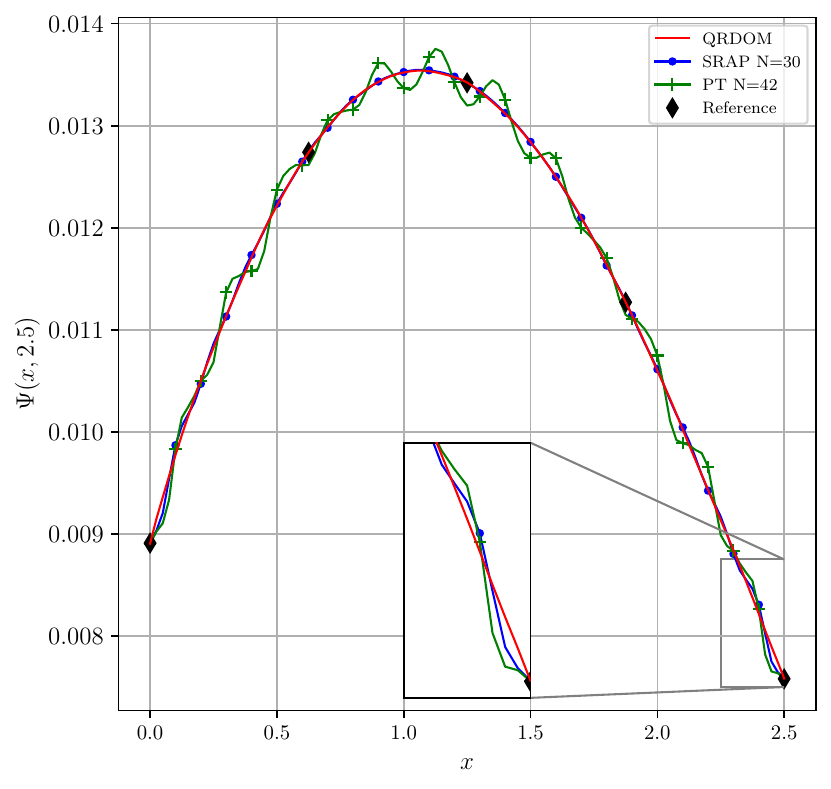}
  \caption{Scalar flux on the top boundary $\Gamma_3$ ($y=2.5$) of the benchmark problem 1.}
  \label{fig:altac-topline}
\end{figure}

The mitigated ray effect solution provided by the QRDOM is notable in
Figure~\ref{fig:altac-topline}, where the profile of the scalar flux
in the top wall ($y=2.5$) is plotted from the solution of the QRDOM,
classical source iteration DOM solutions with the $SRAP_{N}$
\cite{Li1998a} and the $P_NT_N$ \cite{Longoni2001a} quadrature sets,
and the reference values. The reference values were taken from
\cite[Table 4]{Altac2004a}, and both QRDOM and DOM solutions where
computed by a finite element approximation on a uniform mesh of
$128\times 128$ cells. Both methods were initialized with null scalar
flux and used $TOL=10^{-3}$ as stop criteria. The QRDOM convergence
were achieved after $11$ epochs with an average of about $436$ sample
directions on the first octant per epoch. Therefore, the DOM solution
with $SRAP_N$ quadrature set was obtained setting its order to $N=30$,
which gives $495$ discrete directions on the first octant. The DOM
solution with the $P_NT_N$ quadrature set used $441$ discrete
directions on the first octant by setting its order to $N=42$. The highlighted region has a zoom of factor 2.

\subsection{Benchmark problem 2: reflective boundaries}

The benchmark problem 2 has the following parameters: $a=b=1$, $\sigma_s=\sigma_t\equiv 1$, the domain source $Q\equiv 1$ for $x,y\leq 0.52$ and $Q\equiv 0$ otherwise, on the boundaries $Q_{b1}=Q_{b2}=Q_{b3}=Q_{b4}\equiv 0$ and semi-reflective walls $\rho_1=\rho_4\equiv 1$ and $\rho_2=\rho_3\equiv 0$. For this problem the QRDOM has been applied with the scalar flux value at the point $x=y=0.52$ as the goal functional, and $TOL=10^{-4}$. 

Table~\ref{tab:qrdom-loy} shows the QRDOM computed scalar fluxes at
the selected points $(x,y)=(0.5,0.5)$,
$(x,y)=(0.52,0.52)$,$(x,y)=(0.7,0.7)$ and $(x,y)=(0.98,0.98)$. The
obtained values can be compared against the solutions reported in
\cite[Table II]{Loyalka1976a} and the more precise solution from~\cite[Table 5]{Azevedo2018a}.

\begin{table}[h]
  \centering
  \caption{Scalar flux values with $QRDOM$ for problem 2: $a=b=1.0$, $\sigma_s=1.0$. Reference values: $\Psi(0.5, 0.5)=0.687407$, $\Psi(0.7,0.7)=0.344820$ from \cite[Table II]{Loyalka1976a} and $\Psi(0.98,0.98)=0.1326418$ from \cite[Table 5]{Azevedo2018a}. Goal functional $\Psi(0.52,0.52)$ and $TOL = 10^{-5}$.}
  \begin{tabular}[H]{l|ccccc}\hline
    \#cells & $\overline{\Psi}$ & $\Psi(0.5,0.5)$ & $\Psi(0.52,0.52)$ & $\Psi(0.7,0.7)$ & $\Psi(0.98,0.98)$ \\\hline
    $50\times 50$    & $5.9531\cdot 10^{-1}$ & $6.8714\cdot 10^{-1}$ & $6.1281\cdot 10^{-1}$ & $3.4294\cdot 10^{-1}$ & $1.3281\cdot 10^{-1}$ \\
    $100\times 100$  & $5.9536\cdot 10^{-1}$ & $6.8261\cdot 10^{-1}$ & $6.1282\cdot 10^{-1}$ & $3.4304\cdot 10^{-1}$ & $1.3265\cdot 10^{-1}$ \\
    $200 \times 200$ & $5.9537\cdot 10^{-1}$ & $6.8294\cdot 10^{-1}$ & $6.1282\cdot 10^{-1}$ & $3.4304\cdot 10^{-1}$ & $1.3266\cdot 10^{-1}$ \\ \hline
    \end{tabular}
    \label{tab:qrdom-loy}
\end{table}

\begin{figure}[h]
  \centering
  \includegraphics[width=0.7\textwidth]{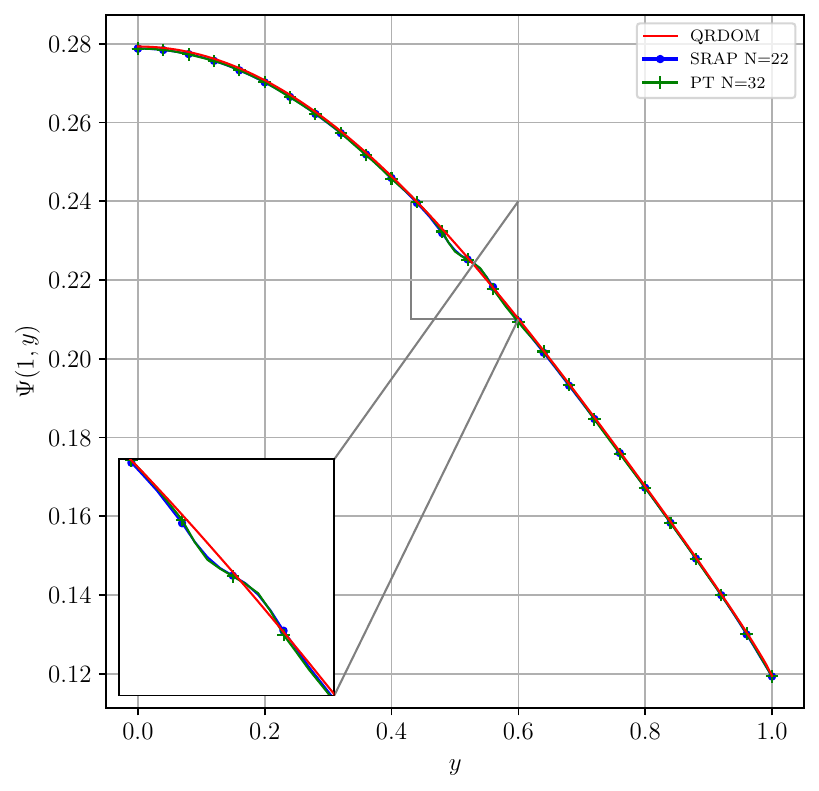}
  \caption{Scalar flux at the right boundary $x=1$ of the benchmark problem 2.}
  \label{fig:loy-rightline}
\end{figure}

It is well known that the classical DOM method applied to this problem will strongly suffer from the ray effect at regions near the lines $x=0.52$ or $y=0.52$. Is is notable that the QRDOM can mitigate this phenomenon as one can observe in Figure \ref{fig:loy-rightline}. This figure shows the profile of the scalar flux at the right wall $x=1$ computed from the $QRDOM$ and from the classical $DOM$ with the $SRAP_N$ and the $P_NT_N$ quadrature sets in a uniform mesh of $100\times 100$ cell and with $TOL=10^{-3}$. To achieve convergence, the QRDOM took $12$ epochs and required an average of $254$ ordinate directions on the first octant per epoch. The classical $DOM$ with both $SRAP_N$ with order $N=22$ and the $P_NT_N$ with order $N=32$ achieved the convergence after $15$ source iterations and used $275$ and $256$ discrete directions on the first octant, respectively.

\subsection{Benchmark problem 3: heterogeneous media}

The third benchmark problem investigated has the following parameters: $a=b=30$, a heterogeneous media with parameters $\sigma_t=1$, $\sigma_s=0.5$, $Q=1$ in the region $x,y\leq 10$, and $\sigma_t=2$, $\sigma_s=0.1$, $Q=0$ otherwise, on the boundaries $\rho_1=\rho_4=1$ and $\rho_2=\rho_3=0$ and null sources $Q_{b1}=Q_{b2}=Q_{b3}=Q_{b4}=0$. The QRDOM solutions were computed by assuming the average scalar flux on the whole domain
\begin{equation}
  \bar{\Psi} = \frac{1}{|\mathcal{D}|}\int_{\mathcal{D}} \Psi\,dx
\end{equation}
as the goal functional and $TOL=10^{-5}$.

\begin{figure}[h]
  \centering
  \includegraphics[width=0.7\textwidth]{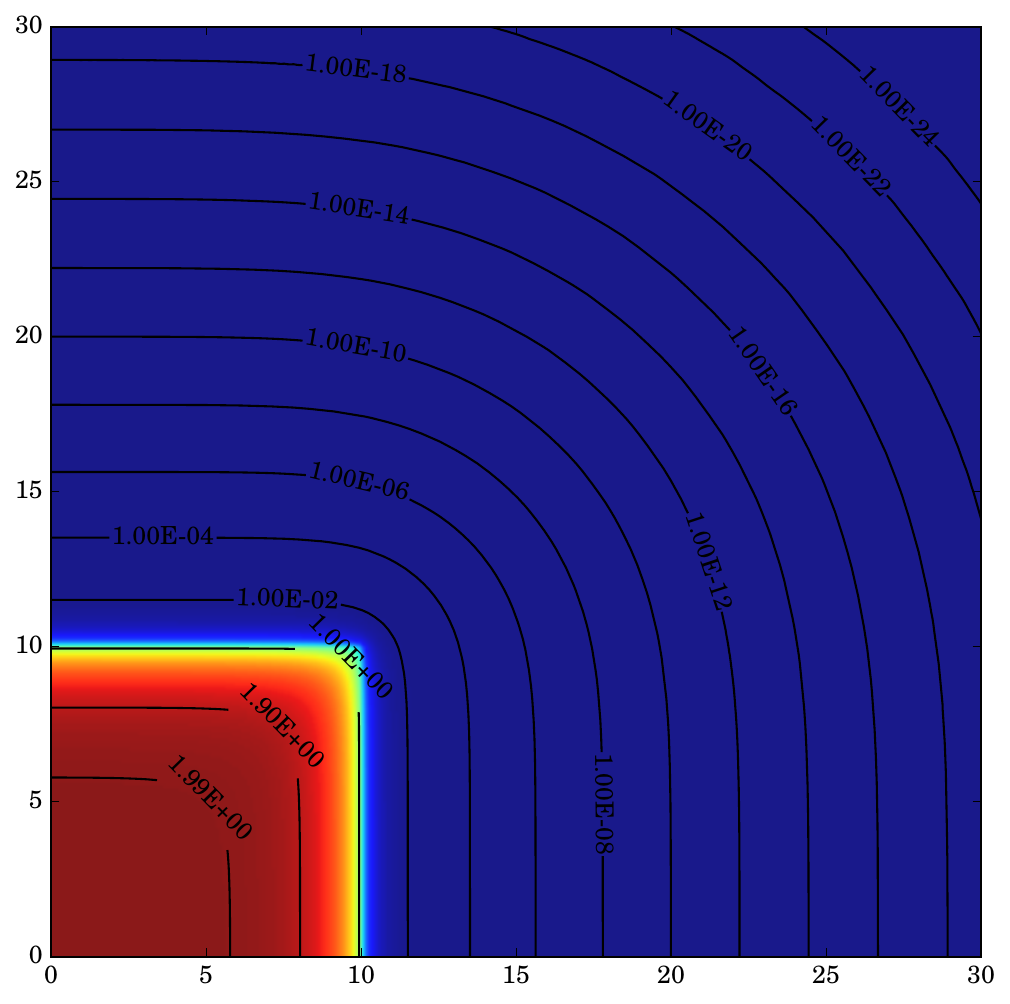}
  \caption{Scalar flux isolines of the benchmark problem 3.}
  \label{fig:bar17a}
\end{figure}

Table \ref{tab:bar17} presents the QRDOM computed values for the goal functional $\bar{\Psi}$ and the average scalar fluxes $\bar{\Psi}_1$ in the region $[0, 10]\times [0, 10]$, $\bar{\Psi}_{2}$ in the region $[10, 30]\times [0, 10]$, $\bar{\Psi}_3$ in the region $[0, 10]\times [10, 30]$, and $\bar{\Psi}_4$ in the region $[10, 30]\times [10, 30]$. Reference values were taken from \cite[Table 8]{Barichello2017a}. In order to highlight the ray effect mitigation, Figure \ref{fig:bar17a} shows the facecolor plot and the isolines of the QRDOM computed scalar flux for this problem.

\begin{table}[h]
  \centering
  \caption{Scalar flux values with $QRDOM$ for problem 3:
    $a=b=30.0$. Reference values \cite[Table 8]{Barichello2017a}:
    $\overline{\Psi}_{1}=1.8360$, $\overline{\Psi}_{2\&3}=1.0678\cdot
    10^{-2}$, $\overline{\Psi}_4=1.1258\cdot 10^{-4}$. Domain mean
    scalar flux as the target functional $\overline{\Psi}$ and $TOL =
    10^{-5}$.} 
  \begin{tabular}[H]{l|ccccc}\hline
    \#cells & $\overline{\Psi}$ & $\overline{\Psi}_{1}$ & $\overline{\Psi}_{2}$ & $\overline{\Psi}_{3}$ & $\overline{\Psi}_4$ \\\hline
    $96\times 96$    &$2.0853\cdot 10^{-1}$ &$ 1.8328$ & $1.0886\cdot 10^{-2}$ & $1.0901\cdot 10^{-2}$ & $1.0782\cdot 10^{-4}$ \\
    $192\times 192$  &$2.0878\cdot 10^{-1}$ &$ 1.8358$ & $1.0695\cdot 10^{-2}$ & $1.0699\cdot 10^{-2}$ & $1.0977\cdot 10^{-4}$ \\
    $384\times 384$  &$2.0883\cdot 10^{-1}$ &$ 1.8364$ & $1.0643\cdot 10^{-2}$ & $1.0660\cdot 10^{-2}$ & $1.1047\cdot 10^{-4}$ \\\hline
    \end{tabular}
    \label{tab:bar17}
\end{table}

\section{Final considerations}\label{sec:final}

In this paper the QRDOM is proposed for the computation of mitigated ray effect approximations of the solutions of one-group transport problems in isotropic mediums. Its central idea is to explore a {\it quasi} Monte Carlo integration within the classical source iteration technique. The application of the QRDOM for problems in rectangular domains was discussed in details and also a convergence acceleration technique has been presented.

The performance of the QRDOM has been tested against three benchmark problems with black walls, reflective walls and in a heterogeneous media. In all cases the method could provide approximations with mitigated ray effects. To this end it demanded hundreds discrete directions on the first octant, which is not sufficient to obtain approximations with mitigated ray effects by the classical DOM with common quadrature sets.

One can observe that the limitations of the presented method are similar to the classical DOM. Its application to transport problems in more complex domains are feasible, including its extension to three-dimensional domains. Moreover, it can be extended to multi-group transport problems as also for problems in anisotropic mediums.

\section*{Acknowledge}

This research has the support of the Centro Nacional de Supercomputação (CESUP) of the Universidade Federal do Rio Grande do Sul (UFRGS).


\end{document}